\documentclass[leqno,12pt]{amsart}
\usepackage{amsfonts}
\usepackage{amsmath,amssymb,amsthm}

\newcommand{\R}{\mathbb R}

\newcommand{\E}{\mathbb E}

\renewcommand{\span}{\mathrm{span}}
\newcommand{\tr}{\mathrm{tr}}

\newtheorem{thm}{Theorem}[section]

\newtheorem{prop}[thm]{Proposition}
\theoremstyle{definition}

\theoremstyle{remark}

\newcommand{\ds}{\displaystyle}

\begin{document}

\title[MARGINALLY TRAPPED MERIDIAN SURFACES OF PARABOLIC TYPE]
{MARGINALLY TRAPPED MERIDIAN SURFACES OF PARABOLIC TYPE IN THE
FOUR-DIMENSIONAL MINKOWSKI SPACE}

\author{Georgi Ganchev and Velichka Milousheva}
\address{Bulgarian Academy of Sciences, Institute of Mathematics and Informatics,
Acad. G. Bonchev Str. bl. 8, 1113 Sofia, Bulgaria}
\email{ganchev@math.bas.bg}
\address{Bulgarian Academy of Sciences, Institute of Mathematics and Informatics,
Acad. G. Bonchev Str. bl. 8, 1113, Sofia, Bulgaria; "L. Karavelov"
Civil Engineering Higher School, 175 Suhodolska Str., 1373 Sofia,
Bulgaria} \email{vmil@math.bas.bg}


\begin{abstract}
A marginally trapped surface in the four-dimensional Minkowski
space is a spacelike surface whose mean curvature vector is
lightlike at each point.  We introduce meridian surfaces of
parabolic type as one-parameter systems of meridians of  a
rotational hypersurface with lightlike axis in Minkowski 4-space
and find their basic  invariants. We find all marginally trapped
meridian surfaces of parabolic type and give a geometric construction of these surfaces.
\end{abstract}

\maketitle

\section{Introduction}

The concept of trapped surfaces was introduced by Roger Penrose in \cite{Pen} and is closely related to  the theory of cosmic
black holes playing an important role in general relativity.
These surfaces were defined in order to study global properties of spacetime. In Physics, a surface in the 4-dimensional spacetime is called marginally trapped if it is closed, embedded, spacelike and
its mean curvature vector is lightlike at each point of the surface.
Recently, marginally trapped surfaces have been studied from a mathematical viewpoint. In the mathematical literature,
it is customary to call a surface in a semi-Riemannian manifold  \emph{marginally trapped}
if its mean curvature vector $H$ is lightlike at each point,
and removing the other hypotheses, i.e. the surface does not need to be closed or embedded.

Classification results in four-dimensional Lorentz space forms were
obtained imposing some extra conditions on the mean curvature
vector, the Gauss curvature or the second fundamental form.
For
example, marginally trapped surfaces with positive relative
nullity in Lorenz space forms were classified by  B.-Y. Chen and
J. Van der Veken in  \cite{Chen-Veken-1}. They also proved the
non-existence of marginally trapped surfaces in Robertson-Walker
spaces with positive relative nullity \cite{Chen-Veken-2} and
classified marginally trapped surfaces with parallel mean
curvature vector in Lorenz space forms \cite{Chen-Veken-3}.

Marginally trapped surfaces in Minkowski 4-space which are
invariant under spacelike rotations were classified by S. Haesen
and  M. Ortega  in \cite{Haesen-Ort-2}.  In \cite{Haesen-Ort-1}  they classified
marginally trapped surfaces in Minkowski 4-space which are
invariant under  boost transformations (hyperbolic rotations).  The
classification of marginally trapped surfaces in Minkowski 4-space
which are invariant under a group of screw rotations (a group of
Lorenz rotations with an invariant lightlike direction) was
obtained in \cite{Haesen-Ort-3}.

Surfaces in the 4-dimensional Minkowski space $\R^4_1$ which are  invariant under spacelike rotations, hyperbolic rotations or screw rotations
are the three types of standard rotational
surfaces with two-dimensional axis  known also as rotational surfaces of elliptic, hyperbolic or parabolic type, respectively.
A rotational surface of elliptic type  is an orbit of
a  regular curve under the action of the orthogonal transformations of $\R^4_1$
which leave a timelike plane point-wise fixed. Similarly, a rotational surface of hyperbolic type
is an orbit of a  regular curve under the action of the orthogonal transformations of $\R^4_1$
which leave a spacelike plane point-wise fixed.
A rotational surface of parabolic type is an an orbit of a  regular curve under the action of the orthogonal transformations of $\R^4_1$
which leave a degenerate plane point-wise fixed.
Some classification results for rotational surfaces of elliptic, hyperbolic or parabolic type with classical extra conditions have been obtained.
A classification of all timelike and spacelike hyperbolic rotational surfaces
with non-zero constant mean curvature in the three-dimensional de Sitter space
$\mathbb{S}^3_1$ is given in \cite{Liu-Liu}
and a classification of the spacelike and timelike Weingarten rotational surfaces of the three types in $\mathbb{S}^3_1$
is found in \cite{Liu-Liu-2}.
In \cite{GM3} we described all Chen spacelike
rotational surfaces of hyperbolic or elliptic type.

\vskip 2mm

In \cite{GM6} we studied  marginally trapped surfaces in the
four-dimensional Minkowski space $\R^4_1$ and developed an
invariant theory of these surfaces  based on the principal lines
generated by the second fundamental form. Using the principal
lines, we introduced a geometrically determined moving frame field
at each point of such a  surface and obtained seven invariant
functions which determine the surface up to a motion in $\R^4_1$.

We applied our theory to a special class of spacelike surfaces
lying on rotational hypersurfaces  with timelike or spacelike
axis. We constructed two-dimensional  surfaces which are
one-parameter systems of meridians of the rotational hypersurface
and called  these surfaces \emph{meridian surfaces}.
The geometric construction of the meridian surfaces is different from the construction of the standard rotational
surfaces with two-dimensional axis. Hence, the class of meridian surfaces is a new source of examples of two-dimensional surfaces in
$\R^4_1$.
We found all
marginally trapped  meridian surfaces lying on rotational
hypersurfaces  with spacelike or timelike axis \cite{GM6}.

In the present paper we continue the study of meridian surfaces  considering a rotational hypersurface with
lightlike  axis in $\R^4_1$ and construct two-dimensional surfaces
which are one-parameter systems of meridians of the rotational
hypersurface. We call these surfaces \emph{meridian surfaces of
parabolic type}.  We calculate their basic  invariants and find all
marginally trapped meridian surfaces of parabolic type.
They are described in  Proposition \ref{P:marginally trapped II class}
and Theorem \ref{T:Marginally trapped - general}.
We give a geometric construction of marginally trapped  meridian surfaces of parabolic type.

Summarizing, we can say that we have described all marginally trapped meridian surfaces of elliptic, hyperbolic and parabolic type.

\section{Preliminaries}

Let  $\R^4_1$ be the Minkowski space endowed with the metric
$\langle , \rangle$ of signature $(3,1)$ and $Oe_1e_2e_3e_4$ be a
fixed orthonormal coordinate system in $\R^4_1$, i.e. $e_1^2 =
e_2^2 = e_3^2 = 1, \, e_4^2 = -1$, giving the orientation of
$\R^4_1$. The standard flat metric is given in local coordinates by
$dx_1^2 + dx_2^2 + dx_3^2 -dx_4^2.$

A surface $M^2$ in $\R^4_1$ is said to be
\emph{spacelike} if $\langle , \rangle$ induces  a Riemannian
metric $g$ on $M^2$. Thus at each point $p$ of a spacelike surface
$M^2$ we have the following decomposition
$$\R^4_1 = T_pM^2 \oplus N_pM^2$$
with the property that the restriction of the metric
$\langle , \rangle$ onto the tangent space $T_pM^2$ is of
signature $(2,0)$, and the restriction of the metric $\langle ,
\rangle$ onto the normal space $N_pM^2$ is of signature $(1,1)$.

A surface $M^2$ in $\R^4_1$ is said to be \emph{timelike} if the
induced metric $g$ on $M^2$ is a metric with index 1, i.e. at each
point $p$ of a timelike surface $M^2$ we have the following
decomposition
$$\R^4_1 = T_pM^2 \oplus N_pM^2$$
with the property that the restriction of the metric $\langle ,
\rangle$ onto the tangent space $T_pM^2$ is of signature $(1,1)$,
and the restriction of the metric $\langle , \rangle$ onto the
normal space $N_pM^2$ is of signature $(2,0)$.

Denote by $\nabla'$ and $\nabla$ the Levi Civita connections on $\R^4_1$ and $M^2$, respectively.
Let $x$ and $y$ denote vector fields tangent to $M$ and let $\xi$ be a normal vector field.
Then the formulas of Gauss and Weingarten give a decomposition of the vector fields $\nabla'_xy$ and
$\nabla'_x \xi$ into a tangent and a normal component:
$$\begin{array}{l}
\vspace{2mm}
\nabla'_xy = \nabla_xy + \sigma(x,y);\\
\vspace{2mm}
\nabla'_x \xi = - A_{\xi} x + D_x \xi,
\end{array}$$
which define the second fundamental tensor $\sigma$, the normal
connection $D$ and the shape operator $A_{\xi}$ with respect to
$\xi$. The mean curvature vector  field $H$ of the surface $M^2$
is defined as $H = \ds{\frac{1}{2}\,  \tr\, \sigma}$.

Let $M^2: z=z(u,v), \,\, (u,v) \in \mathcal{D}$ $(\mathcal{D}
\subset \R^2)$ be a local parametrization on a spacelike surface
in $\R^4_1$. The tangent space at an arbitrary point $p=z(u,v)$ of
$M^2$ is $T_pM^2 = \span \{z_u,z_v\}$, where $\langle z_u,z_u
\rangle
> 0$, $\langle z_v,z_v \rangle > 0$ since $M^2$ is spacelike. We use the standard
denotations $E(u,v)=\langle z_u,z_u \rangle, \; F(u,v)=\langle
z_u,z_v \rangle, \; G(u,v)=\langle z_v,z_v \rangle$ for the
coefficients of the first fundamental form
$$I(\lambda,\mu)= E \lambda^2 + 2F \lambda \mu + G \mu^2,\quad
\lambda, \mu \in \R$$ and we set $W=\sqrt{EG-F^2}$. We choose a
normal frame field $\{n_1, n_2\}$ such that $\langle n_1, n_1
\rangle =1$, $\langle n_2, n_2 \rangle = -1$, and the quadruple
$\{z_u,z_v, n_1, n_2\}$ is positively oriented in $\R^4_1$. Then
we have the following derivative formulas:
$$\begin{array}{l}
\vspace{2mm} \nabla'_{z_u}z_u=z_{uu} = \Gamma_{11}^1 \, z_u +
\Gamma_{11}^2 \, z_v + c_{11}^1\, n_1 - c_{11}^2\, n_2;\\
\vspace{2mm} \nabla'_{z_u}z_v=z_{uv} = \Gamma_{12}^1 \, z_u +
\Gamma_{12}^2 \, z_v + c_{12}^1\, n_1 - c_{12}^2\, n_2;\\
\vspace{2mm} \nabla'_{z_v}z_v=z_{vv} = \Gamma_{22}^1 \, z_u +
\Gamma_{22}^2 \, z_v + c_{22}^1\, n_1 - c_{22}^2\, n_2,\\
\end{array}$$
where $\Gamma_{ij}^k$ are the Christoffel's symbols and the functions
$c_{ij}^k, \,\, i,j,k = 1,2$ are given by
$$\begin{array}{lll}
\vspace{2mm}
c_{11}^1 = \langle z_{uu}, n_1 \rangle; & \qquad   c_{12}^1 = \langle z_{uv},
n_1 \rangle; & \qquad  c_{22}^1 = \langle z_{vv}, n_1 \rangle;\\
\vspace{2mm}
c_{11}^2 = \langle z_{uu}, n_2 \rangle; & \qquad  c_{12}^2 = \langle z_{uv},
n_2 \rangle; & \qquad c_{22}^2 = \langle z_{vv}, n_2 \rangle.
\end{array} $$

Obviously, the surface $M^2$ lies in a 2-plane if and only if
$M^2$ is totally geodesic, i.e. $c_{ij}^k=0, \; i,j,k = 1, 2.$ So,
we assume that at least one of the coefficients $c_{ij}^k$ is not
zero.

The second fundamental form $II$ of the surface $M^2$ at a point
$p \in M^2$ is introduced by the following functions
\begin{equation} \notag
L = \ds{\frac{2}{W}} \left|%
\begin{array}{cc}
\vspace{2mm}
  c_{11}^1 & c_{12}^1 \\
  c_{11}^2 & c_{12}^2 \\
\end{array}%
\right|; \quad
M = \ds{\frac{1}{W}} \left|%
\begin{array}{cc}
\vspace{2mm}
  c_{11}^1 & c_{22}^1 \\
  c_{11}^2 & c_{22}^2 \\
\end{array}%
\right|; \quad
N = \ds{\frac{2}{W}} \left|%
\begin{array}{cc}
\vspace{2mm}
  c_{12}^1 & c_{22}^1 \\
  c_{12}^2 & c_{22}^2 \\
\end{array}%
\right|.
\end{equation}
Let
$X=\lambda z_u+\mu z_v, \,\, (\lambda,\mu)\neq(0,0)$ be a tangent
vector at a point $p \in M^2$. Then
$$II(\lambda,\mu)=L\lambda^2+2M\lambda\mu+N\mu^2, \quad \lambda, \mu \in {\R}.$$
The second fundamental
form $II$ is invariant up to the orientation of the tangent space
or the normal space of the surface.

The condition $L = M = N = 0$  characterizes points at which the
space  $\{\sigma(x,y):  x, y \in T_pM^2\}$  is one-dimensional. We
call such points  \emph{flat points} of the surface. These points
are analogous to flat points in the theory of surfaces in $\R^3$.
In \cite{Lane} and \cite{Little} such points are called inflection
points. The notion of an inflection point is introduced for
2-dimensional surfaces in a 4-dimensional affine space
$\mathbb{A}^4$. E. Lane \cite{Lane} has shown that every point of
a surface is an inflection point if and only if the surface is
developable or lies in a 3-dimensional subspace.

We consider surfaces free of flat points, i.e. $(L, M, N) \neq (0,0,0)$.

The second fundamental form $II$ determines conjugate, asymptotic,
and principal tangents at a point $p$ of $M^2$ in the standard
way.  A line $c: u=u(q), \; v=v(q); \; q\in J \subset \R$ on $M^2$
is said to be an \emph{asymptotic line}, respectively a
\textit{principal line}, if its tangent at any point is
asymptotic, respectively  principal.

The second fundamental form $II$ generates  two invariant functions:
$$k =  \frac{LN - M^2}{EG - F^2}, \qquad
\varkappa = \frac{EN+GL-2FM}{2(EG-F^2)}.$$

 The functions $k$ and $\varkappa$ are invariant
under changes of the parameters of the surface and changes of the
normal frame field \cite{GM5}.
 The sign of $k$ is invariant under congruences
 and the sign of $\varkappa$ is invariant under motions
in $\R^4_1$. However, the sign of $\varkappa$ changes under
symmetries with respect to a hyperplane in $\R^4_1$. It turns out that
the invariant $\varkappa$ is the curvature of the normal
connection of the surface. The number of asymptotic tangents at a
point of $M^2$ is determined by the sign of the invariant $k$.

\section{ Meridian surfaces of elliptic, hyperbolic, and parabolic type in  $\R^4_1$} \label{S:Rotational hypersurfaces}

In \cite{GM4} we constructed a family of surfaces lying on a
standard rotational hypersurface in the four-dimensional Euclidean
space $\R^4$. These surfaces are one-parameter systems of
meridians of the rotational hypersurface, that is why we called
them \emph{meridian surfaces}. We described the meridian surfaces
with constant Gauss curvature, with constant mean curvature, and
with constant invariant $k$.

In the four-dimensional Minkowski space there are three types of
rotational hypersurfaces  -  rotational hypersurfaces with
timelike axis, with spacelike axis, and with lightlike axis. In
\cite{GM6} we used the idea from the Euclidean case to construct
special families of two-dimensional  spacelike surfaces lying
 on rotational hypersurfaces in $\R^4_1$  with
timelike or spacelike axis. The construction was the following.

Let $\{e_1, e_2, e_3, e_4\}$ be the standard orthonormal frame in
$\R^4_1$, i.e.  $e_1^2 =
e_2^2 = e_3^2 = 1, \, e_4^2 = -1$.
First we consider the standard rotational hypersurface with timelike axis.

Let $f = f(u), \,\, g = g(u)$ be smooth functions, defined in an
interval $I \subset \R$, such that $f'\,^2(u) - g'\,^2(u) > 0,
\,\, u \in I$. We assume that $f(u)>0, \,\, u \in I$. The standard
rotational hypersurface $\mathcal{M}'$ in $\R^4_1$, obtained by
the rotation of the meridian curve $m: u \rightarrow (f(u), g(u))$
about the $Oe_4$-axis,  is parameterized as follows:
$$\mathcal{M}': Z(u,w^1,w^2) = f(u)\, \cos w^1 \cos w^2 \,e_1 +  f(u)\, \cos w^1 \sin w^2 \,e_2 + f(u)\, \sin w^1 \,e_3 + g(u) \,e_4.$$

The rotational hypersurface $\mathcal{M}'$ is a two-parameter system of meridians. If $w^1 =
w^1(v)$, $w^2=w^2(v), \,\, v \in J, \, J \subset \R$, we
construct a surface $\mathcal{M}'_m$ lying on $\mathcal{M}'$  in the following way:
\begin{equation}  \notag
\mathcal{M}'_m: z(u,v) = Z(u,w^1(v),w^2(v)), \quad u \in I, \, v \in J.
\end{equation}
Since $\mathcal{M}'_m$ is a
one-parameter system of meridians of $\mathcal{M}'$, we call
$\mathcal{M}'_m$ a \textit{meridian surface on $\mathcal{M}'$}.

In a similar way we consider meridian surfaces lying on the
 rotational hypersurface in $\R^4_1$ with spacelike axis.
Let $f = f(u), \,\, g = g(u)$ be smooth functions, defined in an
interval $I \subset \R$, such that $f'\,^2(u) + g'\,^2(u)
>0$, $f(u)>0, \,\, u \in I$.
The rotational hypersurface $\mathcal{M}''$ in $\R^4_1$, obtained
by the rotation of the meridian curve $m: u \rightarrow (f(u),
g(u))$ about the $Oe_1$-axis  is parameterized as follows:
$$\mathcal{M}'': Z(u,w^1,w^2) = g(u) \,e_1 +  f(u)\, \cosh w^1 \cos w^2 \,e_2 +  f(u)\, \cosh w^1 \sin w^2 \,e_3+ f(u)\, \sinh w^1 \,e_4.$$

If $w^1 = w^1(v), \, w^2=w^2(v), \,\, v \in J, \,J \subset
\R$, we construct a surface $\mathcal{M}''_m$ in $\R^4_1$ in
the following way:
\begin{equation}  \notag
\mathcal{M}''_m: z(u,v) = Z(u,w^1(v),w^2(v)),\quad u \in I, \, v \in J.
\end{equation}
We call $\mathcal{M}''_m$ a
\textit{meridian surface on $\mathcal{M}''$}, since $\mathcal{M}''_m$ is a one-parameter system of
meridians of $\mathcal{M}''$.

In \cite{GM6} we found all  marginally trapped meridian surfaces
lying on the rotational hypersurfaces $\mathcal{M}'$ or
$\mathcal{M}''$. We call the meridian surfaces on $\mathcal{M}'$ or $\mathcal{M}''$
\emph{meridian surfaces of elliptic type} or \emph{meridian surfaces of hyperbolic type},
respectively.

\vskip 2mm Now we shall use the same idea to construct families of
two-dimensional  spacelike surfaces lying
 on a rotational hypersurface in $\R^4_1$  with
lightlike axis.

For convenience we shall use the pseudo-orthonormal base $\{e_1,
e_2, \xi_1, \xi_2 \}$  of $\R^4_1$, defined by $ \ds{\xi_1=
\frac{e_3 + e_4}{\sqrt{2}}},\,\, \ds{\xi_2= \frac{ - e_3 +
e_4}{\sqrt{2}}}$. Note that $\langle\xi_1, \xi_1 \rangle =0$, $\langle
\xi_2, \xi_2 \rangle =0$, $\langle \xi_1, \xi_2
\rangle = -1$. The rotational hypersurface with lightlike axis can
be parameterized by
$$\mathcal{M}''': Z(u,w^1,w^2) =  f(u)\, w^1 \cos w^2 \,e_1 +  f(u)\, w^1 \sin w^2 \,e_2+ \left(f(u) \frac{(w^1)^2}{2} + g(u)\right) \xi_1 + f(u) \,\xi_2,$$
where $f = f(u), \,\, g = g(u)$ are smooth functions, defined in
an interval $I \subset \R$, such that $- f'(u)g'(u) >0$, $f(u)>0,
\,\, u \in I$.

Let $w^1 = w^1(v), \, w^2=w^2(v), \,\, v \in J, \,J \subset \R$
and assume that $(\dot{w}^1)^2 + (\dot{w}^2)^2 \neq 0$.
 We consider the surface $\mathcal{M}'''_m$ in $\R^4_1$ defined by
\begin{equation} \label{E:Eq-1}
\mathcal{M}'''_m: z(u,v) = Z(u,w^1(v),w^2(v)),
\end{equation}
where $u \in I, \, v \in J.$ The surface $\mathcal{M}'''_m$,
defined by \eqref{E:Eq-1}, is a one-parameter system of meridians
of the rotational hypersurface $\mathcal{M}'''$  with lightlike
axis. We shall call $\mathcal{M}'''_m$ a \textit{meridian surface
of parabolic type}.

In the present section we shall find all marginally trapped meridian surfaces of
parabolic type.

Without loss of generality we assume that $w^1 = \varphi(v), \,
w^2=v$. Then the surface $\mathcal{M}'''_m$ is parameterized as
follows:
\begin{equation} \label{E:Eq-2}
\mathcal{M}'''_m: z(u,v) = f(u)\, \varphi(v) \cos v \,e_1 + f(u)\,
\varphi(v) \sin v \,e_2+ \left(f(u) \frac{(\varphi(v))^2}{2} +
g(u)\right)\xi_1 + f(u) \,\xi_2.
\end{equation}

First we shall study the parametric $u$-lines and $v$-lines of the
meridian surface of parabolic type.

Let $v = v_0 = const$ and
denote $c = \varphi(v_0)$, $\alpha = \cos v_0$, $\beta = \sin
v_0$. The parametric  $u$-line $v = v_0 = const$ is given by
\begin{equation} \label{E:Eq-3}
c_u: z(u) = c \alpha f(u)\,e_1 + c \beta f(u)\,e_2+ \left( \frac{c^2}{2} f(u)
 + g(u)\right)\xi_1 + f(u) \,\xi_2.
\end{equation}

Using \eqref{E:Eq-3} we calculate the unit tangent vector field $t_{c_u}$ of $c_u$:
$$t_{c_u}=\frac{1}{\sqrt{-2f' g'}} \left(c \alpha f'\,e_1 + c \beta f'\,e_2+ \left( \frac{c^2}{2} f'
 + g'\right)\xi_1 + f' \,\xi_2\right).$$
We denote by $s$ the arc-length of $c_u$ and calculate the derivative
$$\frac{d t_{c_u}}{d s}=\frac{t_{c_u}'}{s'} = \frac{f' g'' - g' f''}{(-2 f' g')^2}
\left(c \alpha f'\,e_1 + c \beta f'\,e_2+ \left(\frac{c^2}{2} f'
 - g'\right)\xi_1 + f' \,\xi_2 \right).$$
Hence
$\ds{\langle \frac{d t_{c_u}}{d s}, \frac{d t_{c_u}}{d s}  \rangle = \frac{(f' g'' - g' f'')^2}{(-2 f' g')^3}}$.
We set
\begin{equation} \label{E:Eq-4}
n_{c_u}=  \frac{1}{\sqrt{-2f' g'}}
\left(c \alpha f'\,e_1 + c \beta f'\,e_2+ \left(\frac{c^2}{2} f'
 - g'\right)\xi_1 + f' \,\xi_2 \right).
 \end{equation}
Note that $n_{c_u}$ is a timelike vector field, since $\langle n_{c_u}, n_{c_u} \rangle = -1$.  Differentiating \eqref{E:Eq-4} with respect to $s$ we get
$$\frac{d n_{c_u}}{d s}=\frac{n_{c_u}'}{s'} = \frac{f' g'' - g' f''}{(-2 f' g')^2}
\left(c \alpha f'\,e_1 + c \beta f'\,e_2+ \left(\frac{c^2}{2} f'
 + g'\right)\xi_1 + f' \,\xi_2 \right).$$
Thus we obtain the formulas
$$\begin{array}{l}
\vspace{2mm}
\ds{\frac{d t_{c_u}}{d s} = \frac{f' g'' - g' f''}{(-2 f' g')^{\frac{3}{2}}} \,  \,n_{c_u}};\\
\vspace{2mm}
\ds{\frac{d n_{c_u}}{d s} = \frac{f' g'' - g' f''}{(-2 f' g')^{\frac{3}{2}}} \,\,t_{c_u}},
\end{array}$$
which imply that the curvature of $c_u$ is $\ds{ \frac{f' g'' - g' f''}{(-2 f' g')^{\frac{3}{2}}}}$.
For each $v = const$ the
parametric lines $c_u$ are congruent in $\R^4_1$. These curves are the meridians of  $\mathcal{M}'''_m$.
We denote $\varkappa_m(u) = \ds{ \frac{f' g'' - g' f''}{(-2 f' g')^{\frac{3}{2}}}}$.

Now let us consider the parametric  $v$-lines of $\mathcal{M}'''_m$.  Let $u = u_0 = const$ and denote $a =f(u_0)$, $b = g(u_0)$.
The corresponding parametric  $v$-line  is given by
\begin{equation} \label{E:Eq-5}
c_v: z(v) = a \varphi(v) \cos v\,e_1 + a \varphi(v) \sin v\,e_2+ \left( a\frac{\varphi^2(v)}{2} + b\right)\xi_1 + a \,\xi_2.
\end{equation}

Using \eqref{E:Eq-5} we calculate the unit tangent vector field $t_{c_v}$ of $c_v$:
$$t_{c_v}=\frac{1}{\sqrt{\dot{\varphi} ^2 + \varphi^2}} \left((\dot{\varphi} \cos v - \varphi \sin v)\,e_1 + (\dot{\varphi} \sin v + \varphi \cos v)\,e_2+
 \varphi \dot{\varphi} \, \xi_1\right).$$
where $\dot{\varphi}$ denotes the derivative with respect to $v$.
Knowing $t_{c_v}$  we calculate the curvature $\varkappa_{c_v}$ of $c_v$ and obtain that
$\varkappa_{c_v} = \ds{\frac{\varphi \ddot{\varphi} - 2 \dot{\varphi}^2 - \varphi^2 }{ a(\dot{\varphi} ^2 + \varphi^2)^{\frac{3}{2}}}}.$

Let us denote $\overline{\kappa}(v) = \ds{\frac{\varphi \ddot{\varphi} - 2 \dot{\varphi}^2 - \varphi^2 }{(\dot{\varphi} ^2 + \varphi^2)^{\frac{3}{2}}}}$.
Then, for each  $u = u_0 = const$  the curvature of the corresponding parametric  $v$-line is expressed as
$\varkappa_{c_v} = \ds{\frac{1}{a}\,\overline{ \kappa}(v)}$, where $a =f(u_0)$.

\vskip 3mm
Now we shall find the coefficients of the first and the second fundamental forms of the meridian surface of parabolic type.
From \eqref{E:Eq-2} we find the tangent vector fields of $\mathcal{M}'''_m$:
\begin{equation} \label{E:Eq-6}
\begin{array}{l}
\vspace{2mm}
z_u = \ds{f' \varphi \cos v \,e_1 + f'
\varphi \sin v \,e_2+ \left(f' \frac{\varphi^2}{2} +
g'\right)\xi_1 + f' \,\xi_2};\\
\vspace{2mm}
z_v = f( \dot{\varphi} \cos v - \varphi \sin v)\,e_1 + f (
\dot{\varphi} \sin v + \varphi \cos v) \,e_2+ f \varphi \dot{\varphi} \, \xi_1.
\end{array}
\end{equation}
Hence,
the coefficients of the first fundamental form of $\mathcal{M}'''_m$
are $$E = - 2 f'(u) g'(u); \quad F = 0; \quad G = f^2(u) (\dot{\varphi}^2(v) + \varphi^2(v)).$$
The first fundamental form is positive
definite, since $- f' g' >0$. So, $\mathcal{M}'''_m$ is a spacelike surface is $\R^4_1$.

Let us denote
$x =\ds{\frac{z_u}{\sqrt{-2f' g'}}},\,\,y = \ds{\frac{z_v}{f \sqrt{\dot{\varphi}^2 + \varphi^2}}}$.
Then $\{x, y\}$ is an orthonormal tangent frame field of $\mathcal{M}'''_m$. We consider the orthonormal normal frame field, defined by
\begin{equation} \label{E:Eq-7}
\begin{array}{l}
\vspace{2mm}
n_1 = \ds{\frac{1}{\sqrt{\dot{\varphi}^2 + \varphi^2}}  \left((\dot{\varphi} \sin v + \varphi \cos v)\,e_1 + (-\dot{\varphi} \cos v + \varphi \sin v)\,e_2+
 \varphi^2  \, \xi_1\right)};\\
\vspace{2mm}
n_2 = \ds{\sqrt{- \frac{f'}{2g'}}  \left( \varphi \cos v \,e_1 + \varphi \sin v \,e_2 + \frac{f' \varphi^2 - 2g'}{2f'} \, \xi_1 +
 \xi_2\right)}.
\end{array}
\end{equation}

Thus we obtain a frame field $\{x,y, n_1,
n_2\}$ of $\mathcal{M}'_m$, such that $\langle n_1, n_1 \rangle =1$,
$\langle n_2, n_2 \rangle =- 1$, $\langle n_1, n_2 \rangle =0$.

Taking into account
\eqref{E:Eq-6}, we calculate the second partial derivatives of $z(u,v)$:
\begin{equation} \label{E:Eq-8}
\begin{array}{l}
\vspace{2mm}
z_{uu} =\ds{f'' \varphi \cos v \,e_1 + f''
\varphi \sin v \,e_2+ \left(f'' \frac{\varphi^2}{2} +
g''\right)\xi_1 + f'' \,\xi_2};\\
\vspace{2mm}
z_{uv} = \ds{f'(\dot{\varphi} \cos v - \varphi \sin v)\,e_1 + f'
(\dot{\varphi} \sin v  + \varphi \cos v)\,e_2+ f' \varphi \dot{\varphi} \,\xi_1};\\
\vspace{2mm}
z_{vv} =  f\left( (\ddot{\varphi} - \varphi) \cos v - 2 \dot{\varphi} \sin v\right)\,e_1 + f\left(
(\ddot{\varphi}- \varphi) \sin v + 2\dot{\varphi} \cos v\right) \,e_2+ f\left(\dot{\varphi}^2 + \varphi  \ddot{\varphi}\right) \, \xi_1.
\end{array}
\end{equation}

Then equalities \eqref{E:Eq-7} and \eqref{E:Eq-8} imply
\begin{equation} \label{E:Eq-9}
\begin{array}{ll}
\vspace{2mm}
c_{11}^1 = \langle z_{uu}, n_1 \rangle = 0; & \qquad c_{11}^2 = \langle z_{uu}, n_2 \rangle =
\ds{\frac{f''g' - g'' f'}{\sqrt{-2f' g'}}};\\
\vspace{2mm}
c_{12}^1 = \langle z_{uv}, n_1 \rangle = 0; & \qquad  c_{12}^2 = \langle z_{uv}, n_2 \rangle = 0;\\
\vspace{2mm}
c_{22}^1 = \langle z_{vv}, n_1 \rangle = \ds{ f \frac{\varphi \ddot{\varphi} - \varphi^2 - 2 \dot{\varphi}^2}{\sqrt{\dot{\varphi} ^2 + \varphi^2}}};
 &
\qquad
c_{22}^2 = \langle z_{vv}, n_2 \rangle =
\ds{- f  \sqrt{-\frac{f'}{2g'}}\,(\dot{\varphi} ^2 + \varphi^2)}.
\end{array}
\end{equation}

Hence, the coefficients of the second fundamental form are:
\begin{equation} \notag
L = 0; \qquad M = \ds{\frac{f'g'' - g'f''}{2f'g'} \, \frac{\varphi \ddot{\varphi} - \varphi^2 - 2 \dot{\varphi}^2}{\dot{\varphi} ^2 + \varphi^2}};
\qquad N =0.
\end{equation}
Then the invariants $k$ and $\varkappa$ of the meridian surface of parabolic type are expressed as follows:
\begin{equation} \notag
k = - \frac{\kappa_m^2(u) \, \overline{\kappa}^2(v)}{f^2(u)}; \qquad \varkappa = 0.
\end{equation}

\vskip 3mm
The equality $\varkappa = 0$ implies that $\mathcal{M}'''_m$ is a surface with
flat normal connection.

Using \eqref{E:Eq-9} we obtain
\begin{equation}\label{E:Eq-10}
\begin{array}{l}
\vspace{2mm}
\sigma(x,x) =  \qquad \qquad \; - \varkappa_m(u)\,\,n_2; \\
\vspace{2mm}
\sigma(x,y) = 0;  \\
\vspace{2mm}
 \sigma(y,y) = \ds{\frac{\overline{\kappa}(v)}{f(u)}}\, \,n_1 - \frac{1}{f(u)} \sqrt{-\frac{f'(u)}{2g'(u)}}\,\, n_2.
\end{array}
\end{equation}

Taking into account
\eqref{E:Eq-10}, we find the Gauss
curvature $K$ and the mean curvature vector field $H$ of $\mathcal{M}'''_m$:
\begin{equation} \label{E:Eq-11} \notag
K = \ds{- \frac{\kappa_m (u)\, |f'(u)|}{f(u) \sqrt{-2f'(u) g'(u)}}};\\
\end{equation}

\begin{equation} \label{E:Eq-12}
H = \ds{ \frac{\overline{\kappa}(v)}{2f(u)}\,\, n_1 - \frac{1}{2} \left(\kappa_m (u)+ \frac{|f'(u)|}{f  \sqrt{-2f'(u) g'(u)}}\right) \,\, n_2}.
\end{equation}

\vskip 2mm
We can distinguish  two special classes of meridian surfaces of parabolic type.

\vskip 2mm
I. $\overline{\kappa}(v) = 0$. In this case $\mathcal{M}'''_m$ is a surface consisting of flat points, since $L=M=N=0$.
It follows from \eqref{E:Eq-7} that for each  meridian surface of parabolic type the next formulas hold:
$$\begin{array}{l}
\vspace{2mm}
\nabla'_x n_1 = 0;\\
\vspace{2mm}
\nabla'_y n_1 = \ds{\frac{\overline{\kappa}(v)}{f \sqrt{\dot{\varphi}^2 + \varphi^2}}
\left( (\varphi \sin v - \dot{\varphi} \cos v) \,e_1 - (\varphi \cos v + \dot{\varphi} \sin v)  \,e_2 - \varphi  \dot{\varphi} \, \xi_1 \right)}.
\end{array}$$
Having in mind that $\overline{\kappa}(v) = 0$ we get $\nabla'_x n_1 = 0;\,\, \nabla'_y n_1 = 0$, which imply that the normal vector field $n_1$ is constant. Hence,
$\mathcal{M}'''_m$ lies in the hyperplane $\R^3_1$ of $\R^4_1$ orthogonal to $n_1$, i.e. $\mathcal{M}'''_m$ lies in $\R^3_1 =  \span \{x,y,n_2\}$,

In the case $\overline{\kappa}(v) = 0$  the mean curvature vector field is:
\begin{equation} \notag
H = \ds{ - \frac{1}{2} \left(\kappa_m (u)+ \frac{|f'(u)|}{f  \sqrt{-2f'(u) g'(u)}}\right)  n_2}.
\end{equation}

Hence,  $\langle H, H \rangle = 0$ if and only if $H =0$. Consequently, there are no marginally trapped meridian surfaces of parabolic type in the class
$\overline{\kappa}(v) = 0$.

\vskip 2mm
II.  $\varkappa_m(u)  = 0$.
In this case $\mathcal{M}'''_m$ is again  a surface consisting of flat points ($L=M=N=0$).
Since  $\varkappa_m(u)  = 0$, without loss of generality we assume that the meridian curve is determined by $f = u, \,\, g = au +b$,
where $a, b$ are constants, $a <0$. Hence, $\mathcal{M}'''_m$ is a 1-parameter system of straight-lines,
i.e. $\mathcal{M}'''_m$ is a ruled surface, parameterized as follows:
\begin{equation} \label{E:Eq-14}
\mathcal{M}'''_m: z(u,v) = u\, \varphi(v) \cos v \,e_1 + u\,
\varphi(v) \sin v \,e_2+ \left(u \frac{(\varphi(v))^2}{2} +
au +b \right)\xi_1 + u \,\xi_2.
\end{equation}

Let us consider the curve $c: z(v) = \ds{ \varphi(v) \cos v \,e_1 +
\varphi(v) \sin v \,e_2+ \left(\frac{\varphi^2(v)}{2} +
a \right)\xi_1 + \,\xi_2}$. Then the ruled surface $\mathcal{M}'''_m$  is given by
\begin{equation} \label{E:Eq-15}
\mathcal{M}'''_m: z(u,v) = P_0 +  u z(v),
\end{equation}
where $P_0 = b \,\xi_1$ is a fixed point in $\R^4_1$.
The tangent space is spanned by the vector fields $z(v)$ and $\dot{z}(v)$,
and obviously the tangent space is one and the same at the points of each fixed generator of $\mathcal{M}'''_m$.
Hence, $\mathcal{M}'''_m$ is a developable ruled surface in $\R^4_1$. The parametrization \eqref{E:Eq-15} shows that $\mathcal{M}'''_m$
 is a cone in $\R^4_1$ determined by the point $P_0$ and the curve $c: z = z(v)$.

We shall describe the marginally trapped meridian surfaces of parabolic type in the special class $\varkappa_m(u)  = 0$.

\begin{prop} \label{P:marginally trapped II class}
Let $\mathcal{M}'''_m$  be a developable meridian surface of parabolic type, defined by \eqref{E:Eq-14}. Then $\mathcal{M}'''_m$ is marginally trapped
if and only if $\overline{\kappa}^2(v) = \ds{-\frac{1}{2a}}$.
\end{prop}

\noindent
\emph{Proof:} In the case $\varkappa_m(u)  = 0$  the mean curvature vector field is:
\begin{equation} \notag
H = \ds{ \frac{\overline{\kappa}(v)}{2u}\,\, n_1 - \frac{1}{2 u} \sqrt{- \frac{1}{2a}}\, \,n_2}.
\end{equation}
The condition $\langle H, H \rangle =0$ is equivalent to $\overline{\kappa}^2(v) = \ds{-\frac{1}{2a}}$.
\qed

\vskip 3mm
Further we shall consider general meridian surfaces of parabolic type, i.e. we assume that  $\overline{\kappa}(v) \neq 0$ and  $\varkappa_m(u) \neq 0$.

\begin{thm}\label{T:Marginally trapped - general}
The  general meridian surface of parabolic type $\mathcal{M}'''_m$
 is marginally trapped if and only if
$\overline{\kappa}(v) = a = const, \; a \neq 0$, and the
meridian curve is defined by
\begin{equation}  \notag
\begin{array}{l}
\vspace{2mm}
f(u) = u;\\
\vspace{2mm}
g(u) = \ds{\frac{\pm 1}{2a^3}\left(\frac{a^2 u^2 \mp 2 auc}{c \mp au} -2c \ln |c\mp au| + b \right)},
\end{array}
\end{equation}
where $b$ and $c$ are constants, $c \neq 0$.
\end{thm}

\noindent
\emph{Proof:}
Using that
$\varkappa_m(u) = \ds{ \frac{f' g'' - g' f''}{(-2 f' g')^{\frac{3}{2}}}}$ from \eqref{E:Eq-12} we obtain that the mean curvature vector field is
\begin{equation} \notag
H = \ds{ \frac{\overline{\kappa}}{2f}\,\, n_1 + \frac{f(f'' g' - f' g'') + 2 f' g' |f'|}{2f (-2f'g')^{\frac{3}{2}}} \,\, n_2}.
\end{equation}
Hence, $\langle H, H \rangle =0$ if and only if
$$ \overline{\kappa}^2(v) = \frac{\left(f(u)(f''(u) g'(u) - f'(u) g''(u)) + 2 f'(u) g'(u) |f'(u)|\right)^2}{ (-2f'(u)g'(u))^3 }.$$
The last equality implies
\begin{equation}  \label{E:Eq-17}
\begin{array}{l}
\vspace{2mm}
\overline{\kappa}(v) = a = const, \qquad a \neq 0; \\
\vspace{2mm}
\ds{\frac{f(f'' g' - f' g'') + 2 f' g' |f'|}{(-2f'g')^{\frac{3}{2}}}} = \pm a.
\end{array}
\end{equation}

Assuming that the meridian curve
 is given by $f = u; \,\, g = g(u)$,
from  equation \eqref{E:Eq-17} we get
\begin{equation}  \label{E:Eq-17-a}
-u g'' + 2 g' = \pm a (-2 g')^{\frac{3}{2}}.
\end{equation}
After the change
$\ds{\frac{1}{\sqrt{-2 g'(u)}} =  h(u)}$ the above equation is transformed into
\begin{equation}  \label{E:Eq-18}
h' + \frac{h}{u} \pm \frac{a}{u} = 0.
\end{equation}
The general solution of equation \eqref{E:Eq-18} is given by
\begin{equation} \notag
h(u)= \frac{c \mp au}{u}, \qquad c = const, \,\, c \neq 0.
\end{equation}
Hence,
\begin{equation}  \label{E:Eq-19}
g'(u) = \frac{-u^2}{2(c\mp au)^2}.
\end{equation}

Integrating \eqref{E:Eq-19} we obtain that all solutions of
differential equation \eqref{E:Eq-17-a} are given by the formula
\begin{equation}  \label{E:Eq-20} \notag
g(u) = \ds{\frac{\pm 1}{2a^3}\left(\frac{a^2 u^2 \mp 2 auc}{c \mp au} -2c \ln |c\mp au| + b \right)},
\end{equation}
where $b$ and $c$ are constants, $c \neq 0$.
\qed

\vskip 3mm
Theorem \ref{T:Marginally trapped - general} and Proposition \ref{P:marginally trapped II class}
give all marginally trapped meridian surfaces of parabolic type.

\section{Geometric description of marginally trapped meridian surfaces of parabolic type}

In this section we give a complete geometric description of the marginally trapped
meridian surfaces of parabolic type.

\vskip 2mm Note that for the class of marginally trapped meridian
surfaces (in both general and special case) we have
$\overline{\kappa}(v) = a = const, \; a \neq 0$. Now we shall
clear up the geometric meaning of the condition
$\overline{\kappa}(v) = const.$

Each parametric  $v$-line $u = u_0 = const$ of the meridian
surface $\mathcal{M}'''_m$ is given by \eqref{E:Eq-5}. Let us
consider the curve $\overline{c}: \overline{z} = \overline{z}(v)$,
defined by
\begin{equation} \label{E:Eq-21}
\overline{c}: \overline{z}(v) = \varphi(v) \cos v\,e_1 +
\varphi(v) \sin v\,e_2+ \frac{\varphi^2(v)}{2} \,\xi_1 + \xi_2.
\end{equation}
Then each parametric $v$-line is expressed as
\begin{equation} \notag
c_v: z(v) = f(u_0) \,\overline{z}(v) + g(u_0) \,\xi_1.
\end{equation}

Hence, all parametric $v$-lines of $\mathcal{M}'''_m$ are
generated by the curve  $\overline{c}$.

Note that $\langle \overline{z}(v), \overline{z}(v) \rangle =0$.
The curve $\overline{c}$ lies on the paraboloid $\mathcal{P}^2$,
defined by
$$\mathcal{P}^2: z(w^1,w^2) =  w^1 \cos w^2 \,e_1 +   w^1 \sin w^2 \,e_2+ \frac{(w^1)^2}{2} \, \xi_1 + \xi_2.$$

We shall prove that in the case $\overline{\kappa}(v) = a = const,
\; a \neq 0$ the curve $\overline{c}$  is a plane curve on
$\mathcal{P}^2$.

\subsection{Curves on  $\mathcal{P}^2$ with  constant curvature}

Let $\overline{c}$  be the curve on  $\mathcal{P}^2$, given by \eqref{E:Eq-21}.
It follows from \eqref{E:Eq-21} that the unit tangent vector field $\overline{t}(v)$ of $\overline{c}$ is
\begin{equation} \label{E:Eq-22}
\overline{t}(v) = \frac{1}{\sqrt{\dot{\varphi} ^2 + \varphi^2}} \left((\dot{\varphi} \cos v - \varphi \sin v)\,e_1 + (\dot{\varphi} \sin v + \varphi \cos v)\,e_2+
 \varphi \dot{\varphi} \, \xi_1\right).
\end{equation}
Obviously  $\overline{c}$ is a spacelike curve, since  $\langle \overline{t}(v), \overline{t}(v) \rangle =1$.

We denote by $\overline{s}$ the arc-length of $\overline{c}$ and calculate the derivative
$$\frac{d \overline{t}}{d \overline{s}}=\frac{ \dot{\overline{t}}}{\dot{\overline{s}}} = \frac{\varphi \ddot{\varphi} - 2 \dot{\varphi}^2 - \varphi^2 }{(\dot{\varphi} ^2 + \varphi^2)^{2}}
\left((\dot{\varphi} \sin v + \varphi \cos v)\,e_1 + (- \dot{\varphi} \cos v + \varphi \sin v)\,e_2+
 \frac{\dot{\varphi}^4 + \varphi^3 \ddot{\varphi}}{\varphi \ddot{\varphi} - 2 \dot{\varphi}^2 - \varphi^2}  \, \xi_1\right).$$

Hence, the curvature of  $\overline{c}$ is $\overline{\kappa}(v) =
\ds{\frac{\varphi \ddot{\varphi} - 2 \dot{\varphi}^2 - \varphi^2
}{(\dot{\varphi} ^2 + \varphi^2)^{\frac{3}{2}}}}$.

\begin{prop} \label{P:constant curvature-1}
Let $\overline{c}$  be the curve on  $\mathcal{P}^2$, defined by \eqref{E:Eq-21}.
If $\overline{\kappa}(v) = a = const$, $ a \neq 0$, then  $\overline{c}$ is a plane curve.
\end{prop}

\noindent
\emph{Proof:}
We denote
$$ \overline{n}(v) = \frac{1}{\sqrt{\dot{\varphi} ^2 + \varphi^2}} \left((\dot{\varphi} \sin v + \varphi \cos v)\,e_1 + (- \dot{\varphi} \cos v + \varphi \sin v)\,e_2+
 \frac{\dot{\varphi}^4 + \varphi^3 \ddot{\varphi}}{\varphi \ddot{\varphi} - 2 \dot{\varphi}^2 - \varphi^2}  \, \xi_1\right).$$
Then, we have the formula $\ds{\frac{d \overline{t}}{d
\overline{s}}= \overline{\kappa}\, \overline{n}}$. Calculating the
derivative $\dot{\overline{n}}(v)$ we get

\begin{equation} \label{E:Eq-23}
\dot{\overline{n}}(v) = \overline{\kappa}(v)
\left((- \dot{\varphi} \cos v + \varphi \sin v) \,e_1 -  (\dot{\varphi} \sin v + \varphi \cos v)\,e_2+
 \frac{\alpha}{\varphi \ddot{\varphi} - 2 \dot{\varphi}^2 - \varphi^2}  \, \xi_1\right),
 \end{equation}
where $\alpha = \ds{(\dot{\varphi} ^2 + \varphi^2)^{\frac{3}{2}} \frac{d}{d v}
\left(\frac{\dot{\varphi}^4 + \varphi^3 \ddot{\varphi}}{\sqrt{\dot{\varphi} ^2 + \varphi^2}(\varphi \ddot{\varphi} - 2 \dot{\varphi}^2 - \varphi^2)}  \right)}$.

 Let $\overline{\kappa}(v) = a$, $ a \neq 0$. Then
$\varphi \ddot{\varphi} - 2 \dot{\varphi}^2 - \varphi^2 = a(\dot{\varphi} ^2 + \varphi^2)^{\frac{3}{2}}$, which implies
\begin{equation}  \notag
\begin{array}{l}
\vspace{2mm}
\varphi \ddot{\varphi} -  \dot{\varphi}^2 = (\dot{\varphi} ^2 + \varphi^2) \left(1 + a \sqrt{\dot{\varphi} ^2 + \varphi^2}\right);\\
\vspace{2mm}
\varphi (\ddot{\varphi} +  \varphi) = (\dot{\varphi} ^2 + \varphi^2) \left(2 + a \sqrt{\dot{\varphi} ^2 + \varphi^2}\right); \\
\vspace{2mm}
\varphi \dddot{\varphi} -  \dot{\varphi} \ddot{\varphi} = \dot{\varphi} (\ddot{\varphi} + \varphi) \left(2 + 3 a \sqrt{\dot{\varphi} ^2 + \varphi^2}\right).
\end{array}
\end{equation}
Using the last equalities  by straightforward computation we get
$$ \ds{\frac{d}{d v}
\left(\frac{\dot{\varphi}^4 + \varphi^3 \ddot{\varphi}}{(\dot{\varphi} ^2 + \varphi^2)^2}  \right)} = -a^2 \varphi \dot{\varphi},$$
which implies
\begin{equation} \label{E:Eq-24}
\frac{\alpha}{\varphi \ddot{\varphi} - 2 \dot{\varphi}^2 - \varphi^2} = - \varphi \dot{\varphi}.
\end{equation}
Hence, from equalities \eqref{E:Eq-22}, \eqref{E:Eq-23} and \eqref{E:Eq-24}
we obtain the formulas
$$\begin{array}{l}
\vspace{2mm}
\ds{\frac{d \overline{t}}{d \overline{s}}=
\overline{\kappa}\, \overline{n}};\\
\vspace{2mm}
\ds{\frac{d \overline{n}}{d \overline{s}}= -
\overline{\kappa}\, \overline{t}},
\end{array}$$
which imply that the curve $\overline{c}$ is a plane  curve lying in the plane $\span \{\overline{t}, \overline{n}\}$.
\qed

\vskip 3mm
Let us consider the vector fields $T_1 = \ds{\frac{\dot{\varphi} \, \overline{t} + \varphi \, \overline{n}}{\sqrt{\dot{\varphi} ^2 + \varphi^2}}}$
and $T_2 = \ds{\frac{\varphi \, \overline{t} - \dot{\varphi} \, \overline{n}}{\sqrt{\dot{\varphi} ^2 + \varphi^2}}}$.
In the case $\overline{\kappa}(v) = a$, $ a \neq 0$ we calculate that
$$\begin{array}{l}
\vspace{2mm}
T_1 = \ds{\cos v \,e_1 + \sin v \,e_2+ \left( \varphi + \frac{\varphi}{a \sqrt{\dot{\varphi}^2 + \varphi^2}} \right) \, \xi_1};\\
\vspace{2mm}
T_2 = \ds{-\sin v \,e_1 + \cos v \,e_2  - \frac{\dot{\varphi}}{a \sqrt{\dot{\varphi}^2 + \varphi^2}}  \, \xi_1}.
\end{array}$$
$T_1$ and $T_2$ are unit spacelike vector fields such that $\langle T_1, T_2 \rangle =0$ and $\span \{\overline{t}, \overline{n}\} = \span \{T_1, T_2\}$.
Since $\dot{\varphi}^2 + \varphi^2 \neq 0$, the lightlike vector field $\xi_1$ does not lie in the plane $\span \{T_1, T_2\}$.

Each  curve lying on $\mathcal{P}^2$  admits a parametrization of the form $w^1 = \varphi(v), \,\, w^2 = v$ for some smooth function $\varphi$.
From Proposition \ref{P:constant curvature-1} it follows that each curve on $\mathcal{P}^2$  with constant curvature is a plane section of $\mathcal{P}^2$
with a plane which does not contain  $\xi_1$.
Now we shall prove that the converse statement is also true.

\begin{prop} \label{P:constant curvature-2}
Let $c$  be a curve on  $\mathcal{P}^2$, obtained as a plane section with a plane which does not contain  $\xi_1$.
Then  $c$ has constant curvature.
\end{prop}

\noindent
\emph{Proof:}
We shall use the notations $z_1, z_2, \eta_1, \eta_2$ for the coordinate functions of an arbitrary vector field with respect to the base $\{e_1,
e_2, \xi_1, \xi_2 \}$, respectively.  The paraboloid $\mathcal{P}^2$ has the following coordinate parametric equations:
$$ \mathcal{P}^2:
\begin{array}{l}
\vspace{2mm}
z_1 = w^1 \cos w^2;\\
\vspace{2mm}
z_2 = w^1 \sin w^2;\\
\vspace{2mm}
\eta_1 = \ds{\frac{(w^1)^2}{2}};\\
\vspace{2mm}
\eta_2 = 1.
\end{array}
$$
Note that the paraboloid $\mathcal{P}^2$  lies in the hyperplane
of $\E^4_1$, determined by the equation $\eta_2 = 1$. An arbitrary
plane $\pi$ lying in this hyperplane is defined by an equation of
the following form:
$$\pi: A_0 z_1 + B_0 z_2 + C_0 \eta_1 + D_0 = 0,$$
where $A_0, B_0, C_0, D_0$ are constants. Hence, the plane section
of $\mathcal{P}^2$ with $\pi$ is determined by the equation
\begin{equation} \notag
A_0\, w^1 \cos w^2 + B_0 \,w^1 \sin w^2 + C_0\, \frac{(w^1)^2}{2} + D_0 = 0.
\end{equation}
Since we consider plane sections of $\mathcal{P}^2$  with planes
which does not contain  $\xi_1$, we assume that $C_0 \neq 0$ and $A_0^2 + B_0^2 - 2 C_0 D_0 > 0$.
We denote $A = \frac{A_0}{C_0}$, $B = \frac{B_0}{C_0}$, $C = \frac{D_0}{C_0}$ and obtain the equation
\begin{equation} \notag
A \,w^1 \cos w^2 + B \,w^1 \sin w^2 +  \frac{(w^1)^2}{2} + C = 0,
\end{equation}
 or equivalently
\begin{equation} \label{E:Eq-25}
\frac{(w^1)^2}{2} + (A \, \cos w^2 + B \, \sin w^2)\, w^1 + C = 0,
\end{equation}
where $A^2 + B^2 - 2C >0$.
The solution of equation \eqref{E:Eq-25} is
\begin{equation} \notag
w^1 = -(A \, \cos w^2 + B \, \sin w^2) \pm \sqrt{(A \, \cos w^2 + B \, \sin w^2)^2 - 2C}.
\end{equation}

Setting $w^1=\varphi(v), \,\ w^2 = v$ we obtain
\begin{equation} \label{E:Eq-26}
\varphi(v) = -(A \, \cos v + B \, \sin v) \pm \sqrt{(A \, \cos v + B \, \sin v)^2 - 2C}.
\end{equation}

Now we have to prove that the function $\varphi(v)$,  given by formula \eqref{E:Eq-26}, satisfies the condition
$$\ds{\frac{\varphi \ddot{\varphi} - 2 \dot{\varphi}^2 - \varphi^2
}{(\dot{\varphi} ^2 + \varphi^2)^{\frac{3}{2}}}} = const.$$

Let us denote $\theta(v) = A \, \cos v + B \, \sin v$. Then $\varphi = -\theta \pm \sqrt{\theta^2 - 2C}$.
By long but straightforward computation we get
$$\dot{\varphi} ^2 + \varphi^2 = \ds{\frac{(A^2+B^2-2C) \left( \theta \mp \sqrt{\theta^2 - 2C}\right)^2}{\theta^2 - 2C}};$$
$$\varphi \ddot{\varphi} - 2 \dot{\varphi}^2 - \varphi^2 = \ds{\frac{2(A^2+B^2-2C)}{\left(\theta^2 - 2C\right)^{\frac{3}{2}}}
\left((2\theta^2 - C)(\pm \theta - \sqrt{\theta^2 - 2C}) \mp 2C \theta\right)}.$$
Using the last two equalities we calculate that
$$\ds{\frac{\varphi \ddot{\varphi} - 2 \dot{\varphi}^2 - \varphi^2
}{(\dot{\varphi} ^2 + \varphi^2)^{\frac{3}{2}}}} = \mp \frac{1}{\sqrt{A^2+B^2-2C}}.$$

Consequently, the plane section of $\mathcal{P}^2$ with $\pi$ is a curve with constant curvature.

\qed

\subsection{Geometric construction}
Let us consider again the meridian surface of parabolic type $\mathcal{M}'''_m$, defined by \eqref{E:Eq-2}.
Each parametric line $c_u$ ($v = v_0 = const$)  lies in the  plane $\span \{t_{c_u}, n_{c_u}\}$.
Note that
$$\begin{array}{l}
\vspace{2mm}
t_{c_u} - n_{c_u} = \ds{\frac{2 g'}{\sqrt{-2f' g'}}\,\, \xi_1};\\
\vspace{2mm}
t_{c_u} + n_{c_u} = \ds{\frac{2 f'}{\sqrt{-2f' g'}} \left(\varphi(v_0) \cos v_0\,e_1 + \varphi(v_0) \sin v_0\,e_2+ \frac{\varphi^2(v_0)}{2} \,\xi_1 + \xi_2\right)}.
\end{array}$$
Hence, for each $v = v_0 = const$ the $u$-line is a plane curve lying in the plane spanned by the lightlike vector fields $\xi_1$ and
$\overline{z}(v_0)$. So, the meridians of $\mathcal{M}'''_m$ are congruent curves lying in the planes $\span\{ \xi_1, \overline{z}(v)\}$, where
$\overline{z}(v)$ is the position vector of the  curve $\overline{c}$ on the paraboloid
$\mathcal{P}^2$.

Now, taking into account Proposition \ref{P:marginally trapped II
class}, Theorem \ref{T:Marginally trapped - general}, Proposition
\ref{P:constant curvature-1} and Proposition \ref{P:constant
curvature-2}, we obtain a complete description of all  marginally
trapped meridian surfaces of parabolic type. They can be constructed as follows.

\vskip 2mm
I. The general case  $\varkappa_m(u) \neq 0$.

\hskip 5mm $\bullet$
Let $\overline{c}: \overline{z} = \overline{z}(v)$ be a curve on the paraboloid
$\mathcal{P}^2$, obtained by the intersection of $\mathcal{P}^2$ with an arbitrary
plane which does not contain $\xi_1$. The position vector
$\overline{z}(v)$ of $\overline{c}$ is given by \eqref{E:Eq-21}, where
$\varphi$ is determined by \eqref{E:Eq-26}.

\hskip 5mm $\bullet$
The one-parameter system of meridian curves lying in the plane
$\span\{ \xi_1, \overline{z}(v)\}$ and given by
\begin{equation}  \notag
\begin{array}{l}
\vspace{2mm}
f(u) = u;\\
\vspace{2mm}
g(u) = \ds{\frac{\pm 1}{2a^3}\left(\frac{a^2 u^2 \mp 2 auc}{c \mp au}
-2c \ln |c\mp au| + b \right)};
\end{array} \quad b =  const, \; c = const \neq 0
\end{equation}
determines a marginally trapped meridian surface of parabolic type.

\vskip 2mm
II. The special case  $\varkappa_m(u) = 0$.

In this case $\overline{c}$ is again a plane curve on the paraboloid
$\mathcal{P}^2$, obtained by the intersection of $\mathcal{P}^2$ with an arbitrary
plane which does not contain $\xi_1$.
The marginally trapped meridian surface, which is a 1-parameter system of straight-lines, lies in the three-dimensional
space $\R^3_1$ spanned by the plane of the curve $\overline{c}$ and the fixed point $P_0$. Hence, the surface  $\mathcal{M}'''_m$ is
a cone in $\R^3_1$ with lightlike normal vector field.

\end{document}